\documentclass[11pt,reqno]{amsart}

\usepackage{amssymb,latexsym,verbatim,pdfsync,color}

\newcommand{\bbC}{{\mathbb C}}
\newcommand{\bbD}{{\mathbb D}}

\newcommand{\bbR}{{\mathbb R}}

\newcommand{\bbZ}{{\mathbb Z}}

\def\cF{{\mathcal F}}

\def\cH{{\mathcal H}}

\def\cM{{\mathcal M}}

\def\cR{{\mathcal R}}
\def\cS{{\mathcal S}}

\def\cX{{\mathcal X}}

\def\Re{\operatorname{Re}}

\def\la{\langle}
\def\ra{\rangle}

\def\eps{\varepsilon}
\def\z{\zeta} 
\def\vp{\varphi}
\def\ov{\overline}
\def\p{\partial}

\def\ms{\medskip}

\def\MS{M\"untz--Sz\'asz}
\def\tnt{\textstyle{\frac{n}{2}}}

\def\Hol{\operatorname{Hol}}

\newtheorem{thm}{Theorem}[section]
\newtheorem{prop}[thm]{Proposition}
\newtheorem{cor}[thm]{Corollary}

\newtheorem{defn}[thm]{Definition}

\begin{document}

\title[Function
  spaces of exponential growth]{On some spaces of
    holomorphic functions of exponential growth  
on a half-plane} 
\author[M. M. Peloso]{Marco M. Peloso}
\author[M. Salvatori]{Maura Salvatori}
\address{Dipartimento di Matematica, Universit\`a degli Studi di
  Milano, Via C. Saldini 50, 20133 Milano, Italy}
\email{{\tt marco.peloso@unimi.it}}
\email{{\tt maura.salvatori@unimi.it}}
\thanks{{\em Math Subject Classification} 30H99, 46E22, 30C15, 30C40.}
\thanks{Authors partially supported by the grant Prin 2010-11 {\em
    Real and Complex Manifolds: Geometry, Topology and Harmonic
    Analysis}.}

\begin{abstract}
In this paper we study spaces of 
holomorphic  functions   on
the right half-plane $\cR$, 
that we denote by $\cM^p_\omega$,
whose
growth conditions are given in terms of a translation invariant
measure $\omega$ on the closed half-plane $\ov\cR$.
Such a measure has the form $\omega=\nu\otimes m$,  where $m$ is the
Lebesgue measure on $\bbR$ and $\nu$ is a regular Borel measure on
$[0,+\infty)$.   We call these spaces
{\em generalized Hardy--Bergman} spaces on the half-plane $\cR$.

We study in particular the case of 
$\nu$ purely atomic, with point
masses on an arithmetic progression on $[0,+\infty)$. We obtain a
Paley--Wiener theorem for $\cM^2_\omega$, and consequentely the
expression for its reproducing kernel.  We study the growth of
functions in such space and in particular 
show that $\cM^p_\omega$ contains  functions  of order 1.
Moreover, we prove that the orthogonal 
projection from $L^p(\cR,d\omega)$ into $\cM^p_\omega$ is unbounded
for $p\neq2$.

Furthermore, we compare the spaces $\cM^p_\omega$ 
with the classical Hardy and Bergman
spaces, and some other Hardy--Bergman-type  spaces 
introduced more recently.
\end{abstract}
\maketitle

\section{Introduction}

This paper is concerned with spaces of
holomorphic  functions on
the right half-plane $\cR$ whose growth condition is
given in terms of a translation invariant regular measure on $\ov\cR$, and
that can be defined as 
generalized Hardy--Bergman spaces.  
It is easy to see that a measure $\omega$ is  translation invariant on
$\ov\cR$ 
if and only if it has the form $\omega= \nu\otimes m$, where $m$
denotes the Lebesgue measure on $\bbR$ and $\nu$ is a measure on
$[0,+\infty)$. We will simply write $dm(y)=dy$.  These measures play
the same role as the radial measures on the unit disk $\bbD$. \ms

Holomorphic function spaces on $\cR$ with integrability conditions
given in terms of this type of measures have been studied by several
authors, and here we mention in particular Z. Harper \cite{Harper-Laplace,
    Harper-boundedness},  B. Jacob, S. Pott and J. Partington \cite{Jacob-Partington-Pott1},
and I. Chalendar and J. Partington \cite{Chalendar-Partington}, and 
our recent paper \cite{PS}.
We will come back to the spaces these authors considered in Section
\ref{comparison}.\ms

For $0<a<b<\infty$, denote by $S_{a,b}$ the vertical strip
$\{z=x+iy: a<x<b\}$ and by $H^p(S_{a,b})$ the classical Hardy space
$$
H^p(S_{a,b}) =\big\{ f\text{ holomorphic in } S_{a,b}:\,
  \sup_{a<x<b} \int_{-\infty}^{+\infty} |f(x+iy)|^p\, dy
  <\infty\big\}\,.
$$

We simply write $S_b$ to denote $S_{0,b}$.

\begin{defn}\label{mix-spaces-def}{\rm Let $1\le p<\infty$.
Let $\omega=\nu\otimes dy$ be a  translation invariant measure on $\ov\cR$.  

When $\nu(\{0\})=0$, 
we define
\begin{equation}\label{Mp-omega-def1}
\cM^p_\omega(\cR)
= \big\{ f\in \Hol(\cR):\, f\in H^p(S_{a,b}) \text{ for all }
0<a<b<\infty,\ \text{and } f \in L^p(\cR, d\omega)\big\}\, ,
\end{equation}
and if $\nu(\{0\})>0$, 
we define
\begin{equation}\label{Mp-omega-def2}
\cM^p_\omega(\cR)
= \big\{ f\in \Hol(\cR):\, f\in H^p(S_b) \text{ for all }
0<b<\infty,\ \text{and } f \in L^p(\ov\cR, d\omega)\big\} \, .
\end{equation}
In both cases we endow 
$\cM^p_\omega(\cR)$
with the norm
$$
\| f\|_{\cM^p_\omega(\cR)}^p = 
\int_0^{+\infty} \int_\bbR |f(x+iy)|^p\, dyd\nu(x) \, .
$$
We call the spaces $\cM^p_\omega(\cR)$ {\em generalized
  Hardy--Bergman spaces} on the half-plane $\cR$.
}
\end{defn}

We point out that the definition implies that in the case
$\nu(\{0\})>0$, a function in
$\cM^p_{\omega}$, although not initially defined on the imaginary
line, admits a boundary value function that in particular is in
$L^p(dy)$.
\ms

In this paper we consider in particular the measures of the form 
\begin{equation}\label{omega-arho-def}
\omega_{a,\rho} =\nu_{a,\rho}(x) \otimes dy
=\sum_{k=0}^{+\infty} \frac{a^n}{n!} \delta_{\rho\frac{n}{2}}(x) \otimes dy\, ,
\end{equation}
$a,\rho>0$ are fixed parameters.
Thus, the measure 
$\omega_{a,\rho}$ is translation invariant in $\ov\cR$, has
purely atomic part 
$\nu_{a,\rho}$ in the $\Re z$-component and moreover, such atomic measure
$\nu_{a,\rho}$ has support on an arithmetic progression $\{ \frac\rho2
n\}$
with weight $\frac{a^n}{n!}$ at the point-mass $\frac\rho2 n$. 
\ms

Explicitly, $\cM^p_{\omega_{a,\rho}}(\cR)$ 
is the space of holomorphic functions
on $\cR$ that belong to $L^p(\ov\cR,d\omega_{a,\rho}) \cap \big(
\bigcap_{b>0} H^p(S_b))$ with norm
\begin{equation}\label{L2-omega-norm}
\|f\|_{\cM^p_{\omega_{a,\rho}}(\cR)}^p =
\sum_{n=0}^{+\infty} \frac{a^n}{n!} \int_{-\infty}^{+\infty} 
|f\big(\rho\tnt +iy\big)|^p\, dy <\infty \,,
\end{equation}
where  $f$ is defined on the
imaginary axis as its boundary values as function in
$H^p(S_b)$.\ms

We observe that  we equivalently could define
$\cM^p_{\omega_{a,\rho}}(\cR)$ 
as the closure in $L^p(\ov\cR,d\omega_{a,\rho})$ of the space
$\Hol(\ov\cR)\cap \big( \bigcap_{b>0} H^p(S_b))$ with 
norm  given by \eqref{L2-omega-norm}. 
Moreover, it suffices to require that $f\in \bigcap_{n>0}
H^p(S_{\rho\frac{n}{2}})$.\ms  

From simplicity of notation, we write $\cM^p_{a,\rho}$
in place of $\cM^p_{\omega_{a,\rho}}(\cR)$.
\ms

 In \cite{PS} we introduced and studied the space  
$\cM^2_{\omega_{2,1}}=\cM^2_{2,1}$  to give some necessary 
and sufficient conditions for the solutions of the \MS\ problem for
the Bergman space (on the disk $\Delta=D(1,1)$). 
Such problem was first stated and studied by  S. Krantz,
C. Stoppato and the first named author \cite{KPS2}, in
connection with the question of completeness
in the Bergman space on some special domain in $\bbC^2$, the so-called
{\em Diedriech--Forn\ae ss worm domain}.
\ms

 In Section \ref{Muntz-Szasz} we recall the main
results of \cite{PS} (without proofs) to motivated the analysis of the 
slightly more general spaces $\cM^2_{a,\rho}$.  In Section
\ref{M2-arho-sec} 
we prove some of these extensions and, at the same time,
show that some of main properties of $\cM^2_{2,1}$ can
not so easily generalized to the larger class, thus raising some
natural questions.

In Section \ref{comparison} we take a look at other
Hardy--Bergman type spaces studied in particular in  \cite{Harper-Laplace,
    Harper-boundedness}, \cite{Jacob-Partington-Pott1},
and  \cite{Chalendar-Partington} and sometimes  called {\em Zen spaces}.   
These spaces are denoted by
$A^p_\omega$.
We prove a Paley--Wiener theorem for $A^2_\omega$ and, as a
consequence, give a description of the reproducing kernel of $A^2_\omega$.
We also clarify some properties of $A^p_\omega$,
and in particular discuss how functions in such spaces lie in
the Hardy space $H^p(\cR_a)$, for every $a>0$, where $\cR_a$ is the half-plane
$\{\Re z>a\}$.   In this section we also present some examples of
holomorphic functions  and norms on $\cR$ to illustrate some
of the possible behaviors of the average function
$a_{f,p}(x)=\int_\bbR |f(x+iy)|^p\, dy$, for $f\in \Hol(\cR)$.

We conclude the paper with some remarks and open questions.

\ms

\section{The M\"untz--Sz\'asz problem for the Bergman space}
\label{Muntz-Szasz}

\ms

Goal of  this section is to motivate the study of the spaces
$\cM^p_{a,\rho}$ by recalling the main results from \cite{PS},
concerning the special case $\cM^2_{2,1}$.
\ms

Let $\Delta$ be the disk $\{\z:\, |\z-1|<1\}$, denote by $dA$ the
Lebesgue measure in $\bbC$ and consider the (unweighted) Bergman space 
$A^2(\Delta)$.  Then the complex powers $\{\zeta^{\lambda-1}\}$
with $\Re\lambda>0$ are well defined and in $A^2(\Delta)$. 

Following \cite{KPS2}, the \MS\ problem for the Bergman
space is the question of characterizing the sequences $\{\lambda_j\}$ in
$\cR$ such that $\{\z^{\lambda_j-1}\}$ is a complete set in
$A^2(\Delta)$, that is, $\operatorname{span}\{\z^{\lambda_j-1}\}$
is dense in $A^2(\Delta)$.

The classical \MS\ theorem concerns with the completeness of
a set of powers $\{t^{\lambda_j-\frac12}\}$ in $L^2\big([0,1]\big)$, 
where $\Re\lambda_j>0$. 
The solution was provided in two papers separate by M\"untz \cite{Muntz} and by
Sz\'asz \cite{Szasz} where they show that 
the set $\{t^{\lambda_j-\frac12}\}$  is complete $L^2\big([0,1]\big)$ if and only if the
sequence $\{\lambda_j\}$ is a set of uniqueness for the Hardy space 
of the right half-plane $H^2(\cR)$, that is, if $f\in H^2(\cR)$ and
$f(\lambda_j)=0$ for every $j$, then $f$ is identically $0$. 

As in the classical case,
in order to study the \MS\ problem for the Bergman space we 
transformed the question into characterizing the sets of uniqueness for some
(Hilbert) space of holomorphic functions.   In \cite{PS} we showed
that $\{\z^{\lambda_j-1}\}$ is complete in 
$A^2(\Delta)$ if and only if 
$\{\z^{\lambda_j-1}\}$ is a set of uniqueness for
$\cM^2_{2,1}$ and found some sufficient and some necessary
conditions.  We now outline the most relevant results of \cite{PS}.
\ms

For $f\in A^2(\Delta)$ and $z\in\cR$ we
define the Mellin--Bergman transform
\begin{equation}\label{M-Bt}
M_\Delta f(z) = \frac{1}{\pi} \int_\Delta f(\z) \ov \z^{z-1}\, dA(\z)\, .
\end{equation}
The function $\z^{\ov z-1}$ is well defined and belongs to
$A^2(\Delta)$.  Then a set  $\{\z^{\lambda_j-1}\}$ is complete in
$A^2(\Delta)$ if and only if  $f\in A^2(\Delta)$ and
$M_\Delta f(\ov\lambda_j)=0$ for all $j$ implies that $f$ vanishes
identically.  
In order to describe the image of 
$A^2(\Delta)$ under the Mellin--Bergman transform $M_\Delta$
consider the space
\begin{equation}\label{H}
\cH = \big\{ g\in \Hol(\cR):\,  
\textstyle{\frac{\Gamma(1+z)}{2^z}} g(z) \in \cM^2_{2,1} \big\} \, ,
\end{equation}
with norm
\begin{align*}
\|g\|_\cH^2 & = \Big\| \frac{\Gamma(1+z)}{2^z}  g
\Big\|_{\cM^2_{2,1}}^2  
= 
\sum_{n=0}^{+\infty} \int_{-\infty}^{+\infty} |g({\tnt}+iy)|^2
\frac{|\Gamma({\tnt}+1+iy)|^2}{\Gamma(n+1)} \, dy \,. 
\end{align*}

\ms

\begin{thm}\label{new}
The Mellin--Bergman transform
$$
M_\Delta: A^2(\Delta)\to \cH
$$ 
is a surjective isomorphism.
The space $\cH$ consists of holomorphic 
functions on $\cR$ that
are of exponential type at most $\pi/2$ and the polynomials are dense
in $\cH$.  Moreover, it is a Hilbert space with
reproducing kernel
$$
H(z,w) = \frac{1}{2\pi} \frac{\Gamma(z+\ov w)}{\Gamma(1+z)
\Gamma(1+\ov  w)} \, . 
$$
\end{thm}
\ms

We point out that, as corollary of the proof, we obtain a remarkable
factorization theorem for functions in $\cM_{2,1}^2$.  It would be
interesting to prove, if it exists, a similar factorization theorem
for $\cM_{a,\rho}^2$.  

\begin{cor}\label{M2-fact}
We have that
$$
\cM_{2,1}^2 =\frac{\Gamma(1+z)}{2^z} \cH\, .
$$
\end{cor}
\ms

We remark that $\cH$ consists of functions that of exponential type
$\pi/2$ on $\cR$.  Using this factorization we obtained a formula of 
Carleman type for functions in $\cM_{2,1}^2$. 

Recall that the exponent of convergence of a sequence $\{z_j\}$, with 
$|z_j|\to +\infty$, is 
$\rho_1= \inf \{\rho>0:\,  \sum_{j=1}^{+\infty} 1/|z_j|^\rho
<\infty\}$,  the counting function is 
$n(r)  = \# \{ z_j:\, |z_j|\le r\}$ and the upper and lower densities
$d^\pm=d_{\{z_j\}}^\pm$ are 
$$
d^+ =  \limsup_{r\to+\infty} \frac{n(r)}{r^{\rho_1}}\,, \quad
d^- =  \liminf_{r\to+\infty} \frac{n(r)}{r^{\rho_1}}\,.
$$

We are now in the position to state a necessary and a sufficient condition
for zero-sets of $\cM^2_{2,1}$.
\begin{thm}\label{zero-set-thm}
Let $\{z_j\}\subseteq\cR$, $1\le|z_j|\to+\infty$.  The following properties hold.
\begin{itemize}
\item[(i)] 
If $\{z_j\}$ has exponent of convergence $1$
 and  
upper density  $d^+<\frac12$, then
$\{z_j\}$ is a zero-set for $\cM^2_{2,1}\cap\Hol(\ov\cR)$.\smallskip
\item[(ii)]
If $\{z_j\}$ is a zero-set for $\cM^2_{2,1}\cap\Hol(\ov\cR)$, then
\begin{equation}\label{our-Carleman-cond}
\limsup_{R\to+\infty} \frac{1}{\log R} \sum_{|z_j|\le R} \Re\big(1/z_j\big) \le
\frac2\pi  \,.
\end{equation}
\end{itemize}
\end{thm}

The next result gives a partial solution to the 
\MS\ problem for the Bergman space.

\begin{thm}\label{MS-thm}
A sequence $\{z_j\}$ of points in $\cR$ such that $\Re z_j\ge \eps_0$,
for some $\eps_0>0$ and 
that violates condition \eqref{our-Carleman-cond}, is a set of
uniqueness for $\cM^2_{2,1}(\cR)$.

As a consequence, if $\{z_j\}$ is a sequence as above, the set of
powers $\{ \z^{z_j-1}\}$ is a complete set in $A^2(\Delta)$.  
\end{thm}

\ms

\section{The spaces $\cM^2_{a,\rho}$}\label{M2-arho-sec}

\ms

In this section we study the basic properties of $\cM^2_{a,\rho}$.  In
particular we prove a Paley--Wiener type theorem that allows us to
compute its reproducing kernel.  We also prove that the Mellin
transform is a a surjective isometry between a suitable $L^2$ space on
the positive half-line and $\cM^2_{a,\rho}$.

Notice that trivially $H^p(\cR)$ is a closed subset of  $ \cM^p_{a,\rho}$ 
and that $\cM^p_{a,\rho}$
is closed in $L^p(\ov\cR, d\omega)$.

\ms

\subsection{The Paley--Wiener type theorem and its consequences}
\label{PW-sec}
\ms

We being by proving a characterization of $\cM^2_{a,\rho}$ in terms of
the Fourier transform of its boundary values, as in the spirit of the
classical Paley--Wiener theorem.

The Fourier transform of a function $\psi\in L^1(\bbR)$ is
$$
\widehat \psi (\xi) = \frac{1}{\sqrt{2\pi}} \int_{-\infty}^{+\infty}
\psi(x) e^{-ix\xi}\, dx\,. 
$$

For $f\in \cM^2_{a,\rho} $, we 
 write $f_0=f(0+i\cdot)$ to denote its boundary values on the
 imaginary axis. We recall that   the
classical Paley--Wiener theorem for $\cH^2(\cR)$ establishes a
surjective isomorphisms between $\cH^2(\cR)$ and $L^2((-\infty,0),d\xi))$.

\begin{thm}\label{PW-thm} 
Let $f\in \cM^2_{a,\rho}$.  Then $\widehat f_0 \in L^2(\bbR,e^{ae^{\rho\xi}}d\xi)$ and
\begin{equation}\label{PW-cA2-equality}
\|f\|_{\cM^2_{a,\rho}} = \|\widehat f_0
\|_{L^2(\bbR,e^{ae^{\rho\xi}})} \,.
\end{equation}

Conversely, if $\psi\in L^2(\bbR,e^{ae^{\rho\xi}}d\xi)$ and for $z\in\cR$
we set
\begin{equation}\label{PW-cA2-def}
f(z) =\frac{1}{\sqrt{2\pi}} \int_{-\infty}^{+\infty} \psi(\xi)
e^{z\xi}\, d\xi\, ,
\end{equation}
then $f\in \cM^2_{a,\rho}$, equality \eqref{PW-cA2-equality} holds and
$\psi=\widehat f_0$.\ms
\end{thm}

\proof
Since $f\in H^2(S_{\rho\frac{n}{2}})$ for every $n$, by the classical
Paley--Wiener theorem on a strip (\cite{PW}, see also \cite[Thm. 1.1]{PS}), 
we see that 
\begin{align}
\|f\|_{\cM^2_{a,\rho}}^2 
& = \sum_{n=0}^{+\infty} \frac{a^n}{n!} \| 
f(\rho\tnt+i\cdot)\|_{L^2(\bbR)}^2 \notag \\
& = \sum_{n=0}^{+\infty} \frac{a^n}{n!} \big\|
\big[f(\rho\tnt+i\cdot)\big]\widehat{\ }\,\big\|_{L^2(\bbR)}^2 \notag \\ 
& = \sum_{n=0}^{+\infty} \frac{a^n }{n!} \int_{-\infty}^{+\infty}
e^{n\rho\xi}|\widehat f_0(\xi)|^2 \, d\xi \notag \\
& = \int_{-\infty}^{+\infty} 
|\widehat f_0(\xi)|^2 e^{ae^{\rho\xi}}\, d\xi  \label{one-more-PW}
\,.
\end{align}

Conversely, given $\psi\in L^2(\bbR,e^{ae^{\rho\xi}}d\xi)$ 
observe that the integral in \eqref{PW-cA2-def} is absolutely
convergent for
$z\in\cR$:
\begin{equation}\label{pt-eval}
\int_{-\infty}^{+\infty} |\psi(\xi) e^{z\xi}|\, d\xi 
\le \|\psi\|_{L^2(\bbR,e^{ae^{\rho\xi}}d\xi)}  \bigg(\int_{-\infty}^{+\infty}
e^{2x\xi} e^{-ae^{\rho\xi}}\, d\xi \bigg)^{1/2}<\infty\,.
\end{equation}
Therefore, if $f$ is given by \eqref{PW-cA2-def}
it is holomorphic in $\cR$ and 
$f\in H^2(S_{\rho\frac{n}{2}})$ for every $n$, since
$e^{\rho\frac{n}{2}(\cdot)}\psi\in L^2(\bbR)$.
It is also clear that $\widehat f_0=\psi$ and arguing as for 
\eqref{one-more-PW} we obtain \eqref{PW-cA2-equality}.
\qed
\ms

\begin{cor}
The space $\cM^2_{a,\rho}$ is a reproducing kernel Hilbert space and its
reproducing kernel is 
$$
K(z,w) = \frac{1}{2\pi\rho} a^{-\frac{z+\ov w}{\rho}}
\, \Gamma\Big(\frac{z+\ov w}{\rho} \Big)
\,.
$$
\end{cor}

\proof
By \eqref{pt-eval} it follows that point evaluations are continuous in
$\cM^2_{a,\rho}$ and it is elementary to see that it is a Hilbert space.

Let $K_z\in \cM^2_{a,\rho}$ be such that $\la f,\,
K_z\ra_{\cM^2_{a,\rho}}=f(z)$ for every $z\in\cR$ and every $f\in
\cM^2_{a,\rho}$.   
By Theorem \ref{PW-thm} and \eqref{PW-cA2-def} we obtain
\begin{align*}
f(z) 
& = \sum_{n=0}^{+\infty} \frac{a^n}{n!} \big\la f(\rho\tnt+i\cdot),\,
K_z(\rho\tnt+i\cdot)\big\ra_{L^2} \\
& = \sum_{n=0}^{+\infty} \frac{a^n}{n!} \big\la
\big[f(\rho\tnt+i\cdot)\big]\widehat{\ },\,
\big[K_z(\rho\tnt+i\cdot)\big]\widehat{\ }\big\ra_{L^2} \\
& = \sum_{n=0}^{+\infty} \frac{a^n }{n!} \big\la e^{\rho\frac{n}{2}\xi}\widehat f_0,\,
e^{\rho\frac{n}{2}\xi} \widehat K_{z,0} \big\ra_{L^2} \\
& = \int_{-\infty}^{+\infty} \widehat f_0(\xi)
\ov{\widehat {K_{z,0}}(\xi)}\, e^{ae^{\rho\xi}}\, d\xi\,,
\end{align*}
where switching the integral with the sum is justisfied since the last
integral converges absolutely.

On the other hand, 
$$
f(z) =\frac{1}{\sqrt{2\pi}} \int_{-\infty}^{+\infty} e^{z\xi}\widehat f_0(\xi)\, d\xi\,,
$$
so that 
$$
\widehat {K_{z,0}}(\xi) 
= \frac{1}{\sqrt{2\pi}} e^{-ae^{\rho\xi}} e^{\ov z\xi} \,,
$$
and
\begin{align*}
K_z(w) 
&= \frac{1}{2\pi} \int_\bbR e^{w\xi}  e^{-ae^{\rho\xi}}  e^{\ov z\xi}\, d\xi\\
&= \frac{1}{2\pi\rho} a^{-\frac{w+\ov z}{\rho}} \int_0^{+\infty} t^{w+\ov z-1}
e^{-t} \, dt \\
& = \frac{1}{2\pi\rho} a^{-\frac{w+\ov z}{\rho}}
\, \Gamma\Big(\frac{w+\ov z}{\rho} \Big) \,. \qed 
\end{align*}
\ms

It is now clear that $\cM^2_{a,\rho}$ contains functions of order 1 in
the right half-plane, namely, for $w\in\cR$
 $$
K_w(z) = \frac{1}{2\pi\rho} a^{-\frac{z+\ov w}{\rho}}
\, \Gamma\Big(\frac{z+\ov w}{\rho} \Big) \,.
$$
The next result shows that the growth of functions in $\cM^2_{a,\rho}$
is at most of order 1.  Precisely,  
\begin{cor}
Let be $a,\rho>0$ given. 
The functions in 
$\cM^2_{a,\rho}$ satisfy the growth condition
$$
|f(z)| \le C (\Re z)^{1/4} \Big( \frac{2}{a}\Big) ^ {(\Re z)/\rho} \Gamma\Big( \frac{\Re z}{\rho}\Big)\,. 
$$
\end{cor}
\proof
We have
\begin{align*}
|f(z)|
& = \big|
\la \widehat{f_0},\, \widehat{K_{z,0}} \ra_{L^2(\bbR,e^{ae^{\rho\xi}})}
\big|  \le C \| \widehat{K_{z,0}} \|_{L^2(\bbR,e^{ae^{\rho\xi}})} \\
& \le C \bigg( \int_\bbR e^{2
\Re z \xi} e^{-ae^{\rho\xi}} \, d\xi\bigg) ^{1/2}  =C \frac{1}{a^{(\Re z)/\rho}}
  \Gamma \Big( \frac{2\Re z}{\rho}\Big)^{1/2}  
\\ 
& \le C (\Re z)^{1/4} \Big( \frac{2}{a}\Big) ^ {(\Re z)/\rho} \Gamma \Big( \frac{\Re z}{\rho}\Big)
\, , 
\end{align*}
using the standard asymptotics of the Gamma function. 
\qed

\ms

The next result is obvious.

\begin{cor}
 Given $a,a'>0$ and $\rho,\rho'>0$ 
we have the continuous embedding 
$\cM^2_{a,\rho}\subseteq \cM^2_{a',\rho'}$ if and only if $\rho>\rho'$
for any $a,a'$ or if $\rho=\rho'$ and $a> a'$.   
\end{cor}

\ms
There is one more consequence of the Payley--Wiener type theorem about a density 
result in $\cM^2_{a,\rho}$. To state it we first need to introduce
some further notation. 
  
For every $\eps>0$  we 
 denote by $\cM^p_{a,\rho}(\cR_{-\eps})$ the subspace of 
$\cM^p_{a,\rho}$ of functions that are
 holomorphic for $\Re z>-\eps$ and 
that are in $H^p(S_{(-\eps,b)})$ for every $b>0$.  

The following proposition is proved as \cite[Prop. 3.4]{PS} and it is
used in the next section.
\begin{prop}
For $1\le p,q<\infty$, the space
$\bigcap_{\eps>0}\cM^p_{a,\rho}({\cR_{-\eps}})\cap \cM^q_{a,\rho}(\cR_{-\eps})$ is dense in
$\cM^p_{a,\rho}$. 
\end{prop}

\ms

\subsection{The Mellin transform  }
\label{Mellin-subsec}

\ms

We want to show that the Mellin transform is a surjective isometry 
between $\cM^2_{a,\rho}$ and some suitable $L^2$ space of the positive half-line.

More precisely, if $\vp$ is a function defined on $(0,+\infty)$
we consider the (re-normalized)  Mellin transform, that is 
\begin {equation}\label{Mellin-def}
M\vp (z) = \frac{1}{\sqrt{2\pi}} \int_0^{+\infty} \vp(t) t^{z-1}\, dt\,.
\end{equation}

\begin{thm}\label{isometry-thm}
The mapping
$$
M: L^2 \Big( (0,+\infty),\, e^{a\xi^{\rho}}\,
\frac{d\xi}{\xi}\, \Big) \to
\cM^2_{a,\rho}
$$
is a surjective isometry.
\end{thm}

Clearly, the main point of this result is the fact that $M$ is an isometry
that is also surjective.

\proof 
Suppose that $\vp \in L^2 \Big( (0,+\infty),\, e^{a\xi^{\rho}}\,
\frac{d\xi}{\xi}\, \Big) $, then for every $z=x+iy\in\cR$
\begin{align*}
|M\vp(z)|
& \le \frac{1}{\sqrt{2\pi}} \int_0^{+\infty} |\vp(\xi)| \xi^{x-1}\, d\xi\\
& \le \frac{1}{\sqrt{2\pi}} \| \vp \|_{L^2 ( (0,+\infty), e^{a\xi^{\rho}}\,
\frac{d\xi}{\xi})} \Big(
\int_0^{+\infty} e^{-a\xi^{\rho}}  \xi^{2x-1}\, d\xi \Big)^{1/2} \\
& = \frac{1}{\sqrt{2\pi \rho}}  \| \vp \|_{L^2 ( (0,+\infty), e^{a\xi^{\rho}}\,
\frac{d\xi}{\xi})}a^{-\frac{x}{\rho}} \,\Gamma\big(\textstyle{ \frac{2x}{\rho}}\big)^{1/2}
\,. 
\end{align*}
This shows that $M\vp$ is a well defined holomorphic function in $\cR$.
Next, we observe that
\begin{align*}
M\vp(z) 
& = \frac{1}{\sqrt{2\pi}} \int_{-\infty}^{+\infty} \big(\vp\circ\exp)(s) e^{zs}\, ds\\
& = \cF^{-1} \big( (\vp\circ\exp) e^{x(\cdot)}\big)
(y)\,,
\end{align*}
so that we obtain 
\begin{align*}
\big\| M\vp(x+i\cdot) \big\|_{L^2(\bbR)}^2 
& = \| \vp \|_{L^2( (0,+\infty), \xi^{2x-1}d\xi)}^2 \\
& \le C_x  \| \vp \|_{L^2 ( (0,+\infty), e^{a\xi^{\rho}}\,
\frac{d\xi}{\xi})}^2
\, ,
\end{align*}
 uniformly in $x\in (0,b]$.
Hence,  $M\vp \in H^2(S_b)$ for every $b>0$ if $\vp\in L^2 ( (0,+\infty), e^{a\xi^{\rho}}\,
\frac{d\xi}{\xi})$.\ms

Next we consider $\vp,\psi\in  L^2 ( (0,+\infty), e^{a\xi^{\rho}}\,
\frac{d\xi}{\xi})$ and we first assume that both  have compact
support.  Then,  
\begin{align*}
\la M\vp,\, M\psi\ra_{\cM^2_{a,\rho}}
& = \frac{1}{2\pi} \sum_{n=0}^{+\infty} \frac{a^n}{n!} \int_{-\infty}^{+\infty} 
M\vp (\tnt +iy) \ov{M\psi (\tnt +iy) }\, dy \\
& = \sum_{n=0}^{+\infty} \frac{a^n}{n!} \int_{-\infty}^{+\infty} 
\cF^{-1} \Big( \big(\vp\circ\exp) e^{\rho \frac{n}{2}(\cdot)}\Big) (y)
\ov{ \cF^{-1} \Big( \big(\psi\circ\exp) e^{\rho \frac{n}{2}(\cdot)}\Big) (y)} \, dy \\
& = \sum_{n=0}^{+\infty} \frac{a^n}{n!} \int_{-\infty}^{+\infty} 
\vp \big(e^y\big) e^{\rho\frac{n}{2}y} 
\ov{\psi \big(e^y\big) } e^{\rho\frac{n}{2}y}  \, dy \\
& = \int_{-\infty}^{+\infty} 
\vp \big(e^y\big) \ov{\psi \big(e^y\big) } e^{2e^y}\, dy \\
& = \int_0^{+\infty} \vp (t)
\ov{\psi (t) } \, e^{at^{\rho}} \frac{dt}{t}\\
&= \la \vp,\,\psi\ra_{ L^2 ( (0,+\infty), e^{a\xi^{\rho}}\,
\frac{d\xi}{\xi})}
\,. 
\end{align*}
Therefore, the Mellin transform $M$ is a  partial isometry.  
In order to prove that $M$ is onto we need to recall some well-known facts 
about the inversion of the Mellin. 
If, for every $c$ in some interval $I$ it is $g(c+i\cdot)\in L^2(\bbR)\cap L^1(\bbR)$ 
and $g(c+iy)\to 0$ as $|y|\to+\infty$ uniformly in $c\in I$, then 
\begin{equation}\label{M-inv-c}
M_c^{-1} g(\xi) =\frac{1}{\sqrt{2\pi}} 
\int_{-\infty}^{+\infty} g(c+it) \xi^{-c-it} \, dt
\end{equation}
is well defined for every $\xi>0$,  independent of $c\in I$ and $M M_c^{-1}g =g$.
Fix now $\eps>0$ and suppose $f\in \cM^2_{a,\rho}(\cR_{-\eps})\cap \cM^1_{a,\rho}(\cR_{-\eps})$,
satisfying $f(c+iy)\to 0$ as $|y|\to+\infty$ uniformly in $c\in I$.  
Therefore, 
$M_{\rho \frac{n}{2}}^{-1} f$ is independent of $n\ge 0$ and  $M
M_{\rho \frac{n}{2}}^{-1} f=f$.
For every such function $f$, we  define
\begin{equation}\label{def-cM-inv}
M^{-1} f = M_{\frac{n}{2}}^{-1} f\,
\end{equation}
and, by density, it is enough to show that 
\begin{equation}\label{second-isometry}
\|M^{-1} f \|_{L^2( (0,+\infty),\,e^{a\xi^{\rho}}\frac{d\xi}{\xi}\,)} = \| f\|_{\cM^2_{a,\rho}} \,.
\end{equation}

If $f(c+i\cdot)\in  L^2\big(\bbR\big)$ we put $\vp_c(t)=
t^{-c} f\big( c-i\log t)$ and claim that  if $\xi>0$ we have
\begin{equation}\label{M-inv-c-f-equals-M-phi-c}
\xi^c M_c^{-1} f(\xi) =M \vp_c (c+i\log \xi)\,.
\end{equation}
For, 
\begin{align*}
 \big(M \vp_c \big)(c+i\log \xi)
 &=\frac{1}{\sqrt{2\pi}} \int_0^{+\infty} t^{-c} f(c-i\log t)t^{c+i\xi-1}\,dt\\
 &=\frac{1}{\sqrt{2\pi}} \int_{-\infty}^{+\infty} f(c+is) \xi^{-is}\,ds\\
 &=\xi^{c} \big(M_c^{-1} f\big) (\xi)
 \end{align*}
as we claimed. Now,  \cite[Lemma 2.3]{BuJa}  shows that,
for every function $g$ such that $g(t) t^c \in L^2\big((0,+\infty), \frac {dt}t\big)$ it is
$$
\int_0^{+\infty} t^{2c} |g(t)|^2 \, \frac{dt}t= \int_{-\infty}^{+\infty}  |Mg(c+iy)|^2 \,dy\,.
$$
Using all the above, arguing as in \cite[Lemma 3.2]{PS} 
we obtain 
\begin{align*}
\|f\|_{\cM^2_{a,\rho}}^2
& =\| M^{-1} f\|_{L^2( (0,+\infty),\,e^{a\xi^{\rho}}\, \frac{d\xi}{\xi})}^2
\end{align*}
that is, \eqref{second-isometry} holds, and we are done.
\qed

\subsection{Unboundedness of the orthogonal projection}
\ms

Next we consider the question of the  boundedness of the orthogonal
projection.  We find it quite interesting that the orthogonal
projection does not extend to a bounded operator on $\cM_{a,\rho}^p$
for any $p\neq2$.  

\begin{thm}\label{proj-thm}
The orthogonal projection operator 
$P: L^2(\ov\cR,d\omega_{a,\rho} ) \to \cM_{a,\rho}^2$,
 is
unbounded as operator
$$
P: L^p(\ov\cR,d\omega_{a,\rho})\cap L^2(\ov\cR,d\omega_{a,\rho}) \to \cM_{a,\rho}^p
$$
for every $p\neq2$.\ms
\end{thm}

\proof
The proof follows the same lines of the one for the spaces $\cM_{2,1}^p$ 
(see \cite[Thm. 8]{PS}) so that we simply sketch it.
A necessary condition for the projections to be bounded on $L^p$ is that 
the kernels $K_w$ belong to $\cM^{p'}(\ov\cR, d\omega_{a,\rho})$ 
(with $1/p+1/{p'}=1$). By duality and since $P$ is self-adjoint it is enough 
to show that $K_w \notin L^p(\ov\cR,d\omega_{a,\rho})$ for any $p>2$.

First we observe that, there exists $C>0$ such that
$$
\big| \Gamma(\textstyle{\tnt+\frac{u}{\rho}+iy})\big|^p 
 \ge Ce^{-p\frac{u}{ \rho}}  e^{-p\frac{\pi}{2}|y|}  
\exp\Big\{  \textstyle{ \frac{p(n-1)}{2}\log (\frac{n}{2})
  -\frac{pn}{2} } \Big\} \, .
$$

Therefore,
\begin{align*}
\|K_w\|_{\cM^p}^p 
&=C  \frac{1}{a^{\frac{pu}{\rho}} }
\sum_{n=0}^{+\infty} \frac{a^{n(1-\frac{p}{2})}}{n!} \int_{-\infty}^{+\infty} 
\big| \Gamma(\textstyle{\tnt+\frac{u}{\rho}+iy})\big|^p \, dy 
\\
& \ge C_u \sum_{n=0}^{+\infty} \frac{e^
{n[(1-\frac{p}{2}) \log   a-\frac{p}{2}(1+\log2)] }}{n! \, n^{\frac{p}{2}}} e^{\frac{p}{2}n\log n}  \, ,
\end{align*}
which clearly diverges when $p>2$. \

\qed

\ms

\subsection{Some open questions}

We collect here some properties of $\cM^2_{2,1}$ that did not easily
carry over to the more general case of $\cM^2_{a,\rho}$.

\begin{enumerate}
\item Corollary \ref{M2-fact}
shows that
$\cM_{2,1}^2 =\frac{\Gamma(1+z)}{2^z} \cH$, where $\cH$ consists of
functions of exponential type $\pi/2$ in $\cR$.  Thus, functions in 
$\cM_{2,1}^2$ are of order 1 in $\cR$, but can be factored as product
of a non-vanishing term times a function in the space $\cH$.
It is interesting to notice that the space $\cH$, that appears in
\cite{PS} as $M_\Delta(A^2(\Delta))$,  had already appeared
in the literature, in a different context \cite{KT1, KT2}. 

It would certainly be of interest to extend this 
  factorization result to the case $\cM^2_{a,\rho}$, and possibly also
  the description of the analogous of the factor
  $\cH=M_\Delta(A^2(\Delta))$. 

\item  Along the same lines as above, we mention that we do not know
  whether the Carleman formula, Theorem \ref{zero-set-thm} (ii), holds
  in the general case of $\cM^2_{a,\rho}$.   Again, proving  a
factorization theorem
for $\cM_{a,\rho}^2$ would give us a tool for describing its
zero-sets, as we obtained in the case of    $\cM_{2,1}^2$. 

\item It would also be of interest to obtain a description of the \MS\
  problem for {\em weighted} Bergman spaces on $\Delta$.  Here we have
  two most natural choices of weights, namely $v(\z)= (1-|\z-1|^2)^\alpha$, and 
$\tilde v(\z)=|\z|^\alpha$, $\alpha>-1$.  The weight $v$ is radially
symmetric in the disk $\Delta$, while $\tilde v$ is radial in $\bbC$.

\end{enumerate}

\ms

\section{Comparison with other function spaces}\label{comparison} 

\ms

In this section we compare the space $\cM_{a,\rho}^p$ with other
Hilbert spaces of holomorphic functions on the right half-plane.

\subsection{Other spaces of Hardy--Bergman type: Zen
spaces} \label{gen-H-B-spaces-subsec}

Let
\begin{equation}\label{gen-HB-def}
A^p_\omega(\cR)
= \Big\{ f\in \Hol(\cR):  \, \sup_{r>0} \int_0^{+\infty} \int_\bbR
|f(x+r+iy)|^p\, dy\, d\nu(x) <+\infty
\Big\}\, ,
\end{equation} 
where $\omega=\nu(x)\otimes dy$ is a regular Borel measure on $\cR$
and  $\nu$ is  such that there exists $R>0$ such that
\begin{equation}\label{Delta-2-cond}
\sup_{t>0} \frac{\nu([0,2t))}{\nu([0,t))} \le R\, ,
\end{equation}
which is a doubling-type condition at
the origin.

\ms

Such spaces have been studied by several authors; here we mention 
Z. Harper \cite{Harper-Laplace,
    Harper-boundedness},  B. Jacob, S. Pott and J. Partington \cite{Jacob-Partington-Pott1},
and by I. Chalendar and J. Partington \cite{Chalendar-Partington}, and
they are sometimes called {\em Zen spaces}.  We mention here that
previously, this kind of spaces had been considered 
in  \cite{Gustavo} by G. Garrig\'os, 
 in the much more general case of tube domains over cones, although only for
quasi-invariant measures $d\nu$.
\ms

In  \cite{Jacob-Partington-Pott1} two facts about the spaces $A^p_\omega$ are stated
without proof, namely:
\begin{itemize}
\item[(i)] $A^2_\omega$ is a Hilbert space;\smallskip
\item[(ii)]  $f\in A^p_\omega$ implies that $f\in H^p(\cR_a)$ for
  every $a>0$, where $H^p(\cR_a)$ denotes the Hardy space of the
  half-plane $\{\Re z>a\}$.\ms
\end{itemize}

We believe that both statements require some proof.  The reason being 
that, from the definition \eqref{gen-HB-def} it follows that for $f\in
A^p_\omega$, the avarages $a_{f_r,p}=\int_\bbR |f(x+iy)|^p\, dy$ 
are $\nu$-a.e. finite, for every
$r>0$, where
$f_r=f(\cdot+r)$.   
However, this condition is not easily
exploited.

We remark that, for $1\le p<\infty$, the requirement $f\in H^p(S_{a,b})$ for all
$0<a<b<\infty$, implies that the average function $a_{f,p}(x)$
 is finite everywhere
in $(0,+\infty)$ and it is convex.  For, if $x=\delta a+(1-\delta)b$,
for some $0<\delta<1$, then it is well known that
$$
a_{f,p}(x) \le  a_{f,p}(a)^{\delta} a_{f,p}(b)^{1-\delta}
\le \delta  a_{f,p}(a) +(1-\delta) a_{f,p}(b)
 \,.
$$
This implies that
 $a_{f,p}$ must attains its supremum value
as either limit as $x\to 0^+$ or as $x\to+\infty$. 

Hence, once we know that $f\in A^p_\omega$ implies that $f\in H^p(S_{a,b})$,
then the ``sup''-condition in
the definition of the norm of $A^p_\omega$ forces $a_{f,p}(x)\to0^+$
as $x\to+\infty$. Thus, in particular $f\in H^p(\cR_a)$, for every $a>0$.

However, if $f\in\Hol(\cR)$ the average function  $a_{f,p}$ can have
some extreme behaviors, as one could see from
\cite{Harper-Laplace}.  For example, it follows easily from
\cite[Sec. 4]{Harper-Laplace}
that given
any sequence of separated intervals $I_n\subset(0,+\infty)$, that is such that
$\ov I_n\cap \ov I_m
=\emptyset$ for $n\neq m$, there exists $f\in\Hol(\cR)$ such that
$f\in H^p(S_{I_n})$, for all $n$, but $f$ is not in $H^p$ of any
larger strip.

\ms

In the remainder of this section we are going to prove that functions
in $ A^p_\omega$   can be  bounded  uniformly on compact subsets of $\cR$
only in terms of the $ A^p_\omega$-norm and of some constant depending only
the compact set.
  This property  will (easily) imply that 
 $A^2_\omega$ is a Hilbert space, and also 
provide a tool to prove statement (ii) above.  

\begin{prop}\label{sub-mean-val-prop}
Let $1\le p<\infty$ and $ A^p_\omega$ and $\nu$ be as in
\eqref{gen-HB-def} and \eqref{Delta-2-cond}, resp.  Then, for each
compact set $E\subset\cR$ there exists $C_E>0$ such that 
$$
\sup_{\lambda\in E} |F(\lambda)|\le C_E \| F\|_{A^p_\omega} \,.
$$
\end{prop}

\proof
Assume first that the measure $\nu$ has an atomic part $\nu_{\rm at}$.  If
$\nu(\{0\})>0$ then $A^p_\omega$ embeds continuously in $H^p(\cR)$, 
 and the statement follows.  If $\nu(\{0\})=0$, but $\nu_{\rm at}$ has support
 $\{a_n\}_{n\in\bbZ}$, with $a_n\to 0^+$ as $n\to-\infty$.  Arguing as
 before, we obtain that $f\in  A^p_\omega$ belongs to $H^p(\bbC_{a_n})$
 for every $n$, and the conclusion follows again.
\ms

Suppose now that $\nu$ has no atomic part that accumulates at the
orgin.  
Let $a>0$.
By condition \eqref{Delta-2-cond} we can find
$\eps<\frac{a}{4}$  such that 
$\nu([\frac34 \eps, \frac54\eps))>0$.   
Consider the function $g(t) 
=\nu([t-\frac\eps2, t+\frac\eps2))$.  
Notice that $g(\eps)>0$, and
since $\nu$ has no atomic part
that accumulates at the origin, $g$ is continuous, for  $\eps$ 
suitably small.  Then, $g$ has a positive minimum in an interval
$[\eps-\delta,\eps+\delta]$.

For $\lambda \in \ov S_{a-\delta,a+\delta}$, let $\lambda_0 = \lambda
-(a-\eps)$, so that $\Re\lambda_0\in [\eps-\delta,\eps+\delta]$.  
 Denote also $Q(\lambda_0,\eps)$ the square centered in $\lambda_0$ 
with sides of lenght $2\eps$.
   We then
have,
\begin{align}
|F(\lambda)|^p 
& = \frac{1}{\omega(D(\lambda_0,\eps))} \int_{ D(\lambda_0,\eps)}
|F(\lambda_0+a-\eps )|^p\, d\omega(z)\notag \\
& \le \frac{C}{\eps^2 \omega(D(\lambda_0,\eps))} \int_{ D(\lambda_0,\eps)} 
\int_{ D(\lambda_0,\eps)} 
|F(w+a-\eps)|^p\, dA(w) \, d\omega(z)\notag \\
& \le \frac{C}{\eps^2 \omega(D(\lambda_0,\eps))} 
\int_{ D(\lambda_0,\eps)} \int_{ D(z,2\eps)} 
|F(w+a-\eps)|^p\, dA(w) \, d\omega(z)\notag \\
& = \frac{C}{\eps^2 \omega(D(\lambda_0,\eps))} 
\int_{ D(\lambda_0,\eps)} \int_{ D(0,2\eps)} 
|F(z+\z+a-\eps)|^p\, dA(\z) \, d\omega(z)\notag \\
& = \frac{C}{\eps^2 \omega(D(\lambda_0,\eps))} 
 \int_{ D(0,2\eps)} 
\int_{ D(\lambda_0,\eps)} |F(z+\z+a-\eps)|^p\, d\omega(z)\, dA(\z) \notag \\
& \le \frac{C}{\eps^2 \omega(D(\lambda_0,\eps))}
 \int_{ D(0,2\eps)} 
\int_{Q(\lambda_0,\eps)}
|F(x+\xi+a-\eps+i(y+\eta))|^p\, d\nu(x)dy\, dA(\z) \notag \\
& \le C_\eps 
\frac{1}{\nu([\Re\lambda_0-\frac\eps2, \Re\lambda_0+\frac\eps2))}
\|F\|_{A^p_\omega}^p 
\notag \\
& \le C_\eps \,\|F\|_{A^p_\omega}^p  \, .
\label{last-disp-submean}
\end{align}
because of our construction of the function $g$.
\qed

\ms

An easy argument now shows that $A^2_\omega$ is a Hilbert space.  We
do not know of any proof of such fact that does not use the fact that
norm convergence implies uniform convergence on compact subsets.

It also follows that $f\in A^p_\omega$ is bounded on every closed
strip $\ov{S_{a,b}}\subset\cR$.  Therefore, if $f(x+i\cdot)\in
L^p(\bbR)$ for $x=a,b$, it follows that $f\in
H^p(S_{a,b})$.  Hence we have,
\begin{cor}\label{Hp-strip-cor}
{ \rm (1)} If $f\in A^p_\omega$ then $f\in
H^p(\cR_a)$ for any $a>0$ and
$$
\|f\|_{ A^p_\omega}^p
=  \int_0^{+\infty} \int_\bbR |f(x+iy)|^p \, dyd\nu(x)\, .
$$

{\rm (2)} 
When $\nu(\{0\})=0$, we have the equality
$$
A^p_\omega
= \big\{ f\in \Hol(\cR):\, f\in H^p(\cR_a) \text{ for all }
a>0,\ \text{and } f \in L^p(\cR, d\omega)\big\}\, ,
$$
and if $\nu(\{0\})>0$ we have the equality
$$
A^p_\omega
= \big\{ f\in \Hol(\cR):\, f\in H^p(\cR) \text{ and } f \in L^p(\ov\cR, d\omega)\big\} \, .
$$
In both cases we have
$$
\| f\|_{A^p_\omega}^p = 
\int_0^{+\infty} \int_\bbR |f(x+iy)|^p\, dyd\nu(x) \, .
$$
 \end{cor}

\ms

It is now clear that $A^p_\omega$ is a closed subspace of
$\cM^p_\omega$.  It is worth to notice that while functions in
$A^p_\omega$ are bounded in $\cR_a$ for every $a>0$, functions in 
$\cM^p_\omega$ are, in general, allow to grow at infinity, as the case
of $\cM_{a,\rho}^p$ shows.  
\ms

A significant consequence of Prop. \ref{sub-mean-val-prop}
 is a Paley--Wiener theorem type theorem for the space and in
$f\in A^2_\omega$.  In particular it proves that isometry considered
in  \cite[Prop. 2.3]{Jacob-Partington-Pott1}, while 
its proof is inspired by the one of 
\cite[Thm. 2.1]{Harper-Laplace}. 

We need a couple of definitions.  
For $\xi<0$ we set
let
$$
v(\xi) =
\int_0^{+\infty} e^{2\xi x}\, d\nu(x)\, .
$$
It was already observed in \cite[Prop. 2.3]{Jacob-Partington-Pott1}
that the condition \eqref{Delta-2-cond} implies that the integral 
above converges.  For $\vp\in L^2((-\infty,0), v(\xi)d\xi)$ and $z\in\cR$ define
$$
 T\vp (z)
= \frac{1}{\sqrt{2\pi}} \int_{-\infty}^0 e^{z\xi} \vp(\xi)\, d\xi 
= \cF^{-1}\big(e^{x\xi}\vp(\xi)\big)(y)
\, .
$$

\begin{thm}\label{PW-thm-Berg-ab}
{\rm (1)} If $\vp\in  L^2((-\infty,0), v(\xi)d\xi)$ then $T\vp\in
A^2_\omega$ and 
$$
 \|T\vp\|_{A^2_\omega}=  \|\vp\|_{L^2((-\infty,0), v(\xi)d\xi)}
\,. 
$$
{\rm (2)} Conversely, if $F\in A^2_\omega$, then there exists
$\vp\in L^2((-\infty,0), v(\xi)d\xi)$  such that $F=T\vp$ and 
$\| F\|_{A^2_\omega}=
\|\vp\|_{L^2((-\infty,0), v(\xi)d\xi)}$.
\end{thm}

\proof
We remark again that (1) is just
\cite[Prop. 2.3]{Jacob-Partington-Pott1}.

For (2), let $0<a<b<\infty$.
We know that
$\int_\bbR |F(x+iy)|^2\, dy<+\infty$, for all $x>0$. Let $Y>0$ be
fixed and let $R_Y$ be the rectangle of vertices $(a,-Y)$, $(b,-Y)$,
$(b,Y)$,  and $(a,Y)$.  By Cauchy's integral formula, for $\lambda \in
R_Y$ we have
$$
F(\lambda)= \frac{1}{2\pi i} \int_{\p R_Y} \frac{F(z)}{z-\lambda}\,
dz\, .
$$
 Letting $Y\to+\infty$, using the fact that $F$ is bounded in every
 closed strip,  for $a< \Re\lambda< b$  we obtain
\begin{align*}
F(\lambda)
& = \frac{1}{2\pi i} \Big(
\int_{b+i\bbR} \frac{F(z)}{z-\lambda}\, dz - \int_{a+i\bbR} \frac{F(z)}{z-\lambda}\, dz
\Big) \\ 
& = \frac{1}{2\pi } \Big(
\int_\bbR \frac{F(b+iy)}{iy-(\lambda-b)}\, dy 
-\int_\bbR \frac{F(a+iy)}{iy-(\lambda-a)}\, dy \Big)\, .
\end{align*}
Now notice that
$$
\big| \int_\bbR \frac{F(b+iy)}{iy-(\lambda-b)}\, dy \big|^2
\le \| F_b\|_{L^2(\bbR)}^2 \int_\bbR \frac{1}{y^2+|\Re\lambda-b|^2}\, dy
\to +\infty 
$$
as $b\to+\infty$, since $F\in H^2(\cR_{-\eps})$ for every $\eps>0$. 
Therefore,  observing that $\lambda-a\in \cR$,
\begin{equation}\label{who}
F(\lambda)= 
\frac{1}{2\pi } \int_\bbR \frac{F(a+iy)}{(\lambda-a)-iy}\, dy =
\cS(F_{a})(\lambda-a)\, ,
\end{equation}
where $\cS$ denotes the Szeg\"o projection on $H^2(\cR)$.
By the classical Paley--Wiener theorem there
exists $g_a \in L^2((-\infty,0),d\xi)$ such that
$$
F(\lambda)= \cS(F_{a})(\lambda-a)
= \frac{1}{\sqrt{2\pi}} \int_{-\infty}^0 e^{(\lambda-a)\xi} g_a(\xi)\, d\xi\, ,
$$
for $\Re\lambda>a$.  By the
uniqueness of the Fourier transform,
$e^{-a\xi} g_a = e^{-a'\xi} g_{a'}$ for every $a,a'>0$.  

Hence,
\begin{equation}\label{F-from-vp}
F(\lambda)
= \frac{1}{\sqrt{2\pi}} \int_{-\infty}^0 e^{\lambda\xi} \vp(\xi)\, d\xi 
\, ,
\end{equation}
where, for $\xi<0$,
$\vp(\xi) =
 e^{-a\xi} g_a(\xi)$, any $a>0$, so that 
$e^{a\xi}\vp \in L^2((-\infty,0),d\xi)$ for every $a>0$. Thus
 the integral in
\eqref{F-from-vp} is indeed absolutely convergent and
 by Plancherel theorem again,
\begin{align*}
\| F\|_{A^2_\omega}^2
& = \sup_{r>0}  \int_0^{+\infty} \int_\bbR |\widehat{F_{x+r}}(\xi)|^2\, d\xi\, d\nu(x)\\
& = \sup_{r>0}   \int_0^{+\infty} \int_{-\infty}^0 e^{2(x+r)\xi}
|\vp(\xi)|^2\,  d\xi \,
d\nu(x)\\
& = \sup_{r>0}   \int_{-\infty}^0  \Big( \int_0^{+\infty}
e^{2x\xi}\, d\nu(x)\Big)   e^{2r\xi} |\vp(\xi)|^2 
 \, d\xi \\
& =  \int_{-\infty}^0  |\vp(\xi)|^2  v(\xi)\, d\xi\, ,
\end{align*}
as we wished to show.
\ms
\qed

\ms
As usual, given a Paley--Wiener theorem for $A^2_\omega$, it is possible to obtain 
its reproducing kernel. Indeed,
\begin{thm}\label{RK-thm-Berg-ab}
The reproducing kernel for $A^2_\omega$ is 
 $$ 
 K(z,w)=\frac 1 {2\pi} \int_{-\infty}^0 e^{(z+\ov w)\xi} \frac {d\xi}{v(\xi)}.
 $$
 \end{thm}
 
 \proof
 With the notation as in  Theorem \ref{PW-thm-Berg-ab}, to every function $G$ in 
 $A^2_\omega$ we associate the function $\vp_G$ in  $L^2((-\infty,0), v(\xi)d\xi)$
 such that $T\vp_G=G$.

From Proposition \ref{sub-mean-val-prop} we know that  $A^2_\omega$ is
a reproducing kernel Hilbert space.
 Let $K_z$ be the reproducing kernel in $A^2_\omega$, that is, for 
 $F\in A^2_\omega$ and $z\in \cR$ it holds that 
 $F(z)=\la F, K_z\ra_{A^2_\omega}$.

 Then,
 \begin{align*}
F(z) 
& = \int_0^{+\infty} \la F(u+i\cdot), K_z(u+i\cdot)\ra_{L^2(\bbR)} d\nu(u)\\
& = \int_0^{+\infty}  \big\la \big[F(u+i\cdot)\big]\widehat{\ },\,
\big[K_z(u+i\cdot)\big]\widehat{\ }\big\ra_{L^2(\bbR)} d\nu(u)\\
& = \int_0^{+\infty}  \big\la e^{u\xi}\vp_F,\,
e^{u\xi}  \vp_{K_{z} }\big\ra_{L^2(\bbR)} d\nu(u)\\
& = \int_{-\infty}^0 \vp_F(\xi) \ov{ \vp_{K_{z}} (\xi) }v(\xi)\,d\xi,
\end{align*}
where switching the integral with the sum is justisfied since the last
integral converges absolutely.

On the other hand, 
$$
F(z) =\frac{1}{\sqrt{2\pi}} \int_{-\infty}^0 e^{z\xi}  \vp_F(\xi)\, d\xi\,,
$$
so that 
$$
\vp_{K_{z} }(\xi) 
= \frac{1}{\sqrt{2\pi}} \frac{1}{v(\xi)} e^{\ov z\xi} \,,
$$
and
$$
K_z(w) 
= \frac{1}{2\pi} \int_ {-\infty}^0       e^{w\xi}    e^{\ov z\xi}  \frac{1}{v(\xi)}\, d\xi$$
as we wished to show. \qed

\subsection{The space $L^2(\ov\cR,d\omega)\cap \Hol(\ov\cR)$}

We take the opportunity to discuss the naive definition of Bergman
space with respect to a general measure such as $\omega_{a,\rho}$ .  Define
$\cX_{a,\rho}^2$ as the closure of the holomorphic functions on $\cR$ that
extend continuously to the boundary and in the norm of
$L^2(\ov\cR,d\omega_{a,\rho})$. 

We show that such definition is quite inadequate from the point of
view of complex
analysis.  For, we prove the following

\begin{prop}\label{unb-eval-pts}
The following properties hold:
\begin{itemize}
\item[(1)] 
 the set of points $z\in\cR$ such that 
$\cX_{a,\rho}^2\ni f\mapsto f(z)$ is unbounded is dense in $\ov\cR$;\smallskip
\item[(2)]
$\cX_{a,\rho}^2$ contains functions that are not holomorphic
  in $\cR$. \ms
\end{itemize}
\end{prop}
\proof
We prove the statement for the space $\cX_{2,1}^2$, the proof
in the general case being completly analogous.

(1) Let $h(z)= (1+z)^{-1} \exp\{i e^{2\pi iz}\}$  and set
  $f_k(z)=h(kz)$. Observe that
$|\exp\{i e^{2\pi ik(\frac n2 +iy)}\}|= 
|\exp\{ ie^{-2k\pi y}\cos (kn\pi) \}|=1$.  
Then 
\begin{itemize}
\item[{\tiny $\bullet$}] 
$f_k(\tnt+iy)\to 0$ a.e. $dy$  as $k\to +\infty$;\smallskip
\item[{\tiny $\bullet$}]  
$\| f_k\|_{L^2(\cR,d\omega_{2,1})}^2 
= \pi \sum_{n=0}^{+\infty} \frac{2^n}{n!} \big(1+k\frac n2\big)^{-2} \to 0$, as $k\to +\infty$,\ms
\end{itemize}
as an easy calculation shows.  Now let $z\in\cR$ with $\Re z>0$
rational, equal to $p/q$ with $p,q$ relatively prime and $q\neq2$.  
Such points are dense in $\ov\cR$ and
$$
f_k (z)= \frac{1}{1+\frac{kp}{q} +iky} \exp \{ ie^{-2k\pi y} \big(
\cos(kp\pi/q) +i \sin(kp\pi/q) \big)\}\, .
$$
If we choose $k=\ell_0 +2q\ell$, $\ell=1,2,\dots$ we see that
$$
|f_k(z)| \approx \ell^{-1} e^{-2(\ell_0+2q\ell)y}\sin (\ell_0 p\pi/q)\, .
$$
Choosing $\ell_0$ such that $\sin (\ell_0 p\pi/q)$ we see that
$|f_k(z)| \to +\infty$ if $y<0$.  Thus, the point evalution
at $z=\frac pq+iy$ with $y<0$ are not bounded on $\cX_{2,1}^2$.  

To deal with the case $y>0$ it suffices to replace 
$\exp\{i e^{2\pi iz}\}$ with $\exp\{i e^{-2\pi iz}\}$ in the
definition of $f_k$.\ms

(2) Let 
$$
f_k(z) =\frac{1}{1+z} \frac{\exp\{ ie^{4k\pi iz}-1\} }{e^{4k\pi iz}} \,.
$$
Since 
$$
f_k({\tnt}+iy) 
= \frac{1}{1+\frac n2 +iy} \frac{\exp\{ ie^{-4k\pi y}-1\} }{e^{-4k\pi y}}
$$
and 
for $t>0$, $\big| \frac{e^{it}-t}{t}\big|\le C$, we have
\begin{align*}
\| f_k\|_{L^2(\cR,d\omega{2,1})}^2 
& =\sum_{n=0}^{+\infty} \frac{2^n}{n!} \int_\bbR |f_k(\tnt+iy)|^2\, dy \\
& \le C \sum_{n=0}^{+\infty} \frac{2^n}{n!} \int_\bbR
\frac{1}{\big(1+\frac n2\big)^2 +y^2} \, dy\\
& \le C
\,.
\end{align*}
Then, $f_k\in \cH$ for $k=1,2,\dots$.  Moreover,
$f_k\to g$ $\omega$-a.e. where,
$$
g(\tnt+iy)= \begin{cases} 
i(1+\tnt+iy)^{-1} \quad & y>0 \cr
0 & y<0\, .
\end{cases}
$$
It is now easy to see that $f_k\to g$ in $L^2(\bbR,d\omega)$.  Since
$g$ cannot be extended to a holomorphic function on $\cR$, this
concludes the proof.
\qed

\ms

\section{Final remarks}
\ms

It is interesting to notice that the space $\cH$ defined in \eqref{H},
and that in \cite{PS} we showed to be equal to $M_\Delta(A^2\Delta)$, 
had already appeared
in the literature, in a different context \cite{KT1, KT2}.  However,
$\cH$ can be described as the closure of polynomials in the $L^2(\cR,
d\mu)$-norm 
where $d\mu(z)=2^{-x}|\Gamma(1+z)|^2 dA(z)$, which is not translation 
invariant in $\cR$.  In \cite{KT1} the authors also discussed a
\MS-type question, concerning the completeness of the powers
$\{(1-z)^{\lambda_n}\}$ in $H^2(\bbD)$, for $\lambda_n>0$ and
$\lambda_{n+1}-\lambda_n>\delta$; their results however have no
(obvious) connection with the \MS\ problem for the Bergman space.

Of course, there exist other papers dealing with \MS-type questions.  
In \cite{Sed} A. Sedletkskii studied the completeness of sets of
exponentials in weighted $L^p$ spaces on $(0,+\infty)$
in terms of zeros of functions in the classical Bergman space on a
half-plane.   \ms

In \cite{Lukacs57,Lukacs} E. Lukacs studied positive measures $\nu_L$ 
on the real line
that are Fourier transform of restriction of entire functions.  Then, 
we can consider
the measures on $\ov\cR$ of the form $\nu_L\otimes dy$, (with
$\nu_L$ restricted to $[0,+\infty)$).
We believe these
measures constitute an
interesting class of measures for which studying the properties of the
function spaces  $\cM^p_\omega$.  We wish to come back to this, and
the other open problems we mentioned, in a future work.

\ms

\bibliography{m2-mathscinet-1}

\begin{thebibliography}{10}

\bibitem{BuJa}
P.~L. Butzer and S.~Jansche.
\newblock A self-contained approach to {M}ellin transform analysis for square
  integrable functions; applications.
\newblock {\em Integral Transform. Spec. Funct.}, 8(3-4):175--198, 1999.

\bibitem{Chalendar-Partington}
I.~Chalendar and J.~R. Partington.
\newblock Norm estimates for weighted composition operators on spaces of
  holomorphic functions.
\newblock {\em Complex Anal. Oper. Theory}, 8(5):1087--1095, 2014.

\bibitem{Gustavo}
G.~Garrig{\'o}s.
\newblock Generalized {H}ardy spaces on tube domains over cones.
\newblock {\em Colloq. Math.}, 90(2):213--251, 2001.

\bibitem{Harper-boundedness}
Z.~Harper.
\newblock Boundedness of convolution operators and input-output maps between
  weighted spaces.
\newblock {\em Complex Anal. Oper. Theory}, 3(1):113--146, 2009.

\bibitem{Harper-Laplace}
Z.~Harper.
\newblock Laplace transform representations and {P}aley-{W}iener theorems for
  functions on vertical strips.
\newblock {\em Doc. Math.}, 15:235--254, 2010.

\bibitem{Jacob-Partington-Pott1}
B.~Jacob, J.~R. Partington, and S.~Pott.
\newblock On {L}aplace-{C}arleson embedding theorems.
\newblock {\em J. Funct. Anal.}, 264(3):783--814, 2013.

\bibitem{KPS2}
S.~G. Krantz, M.~M. Peloso, and C.~Stoppato.
\newblock Completeness on the worm domain and the {M}\"untz--{S}z\'asz problem
  for the {B}ergman space.
\newblock {\em preprint}, 2015.

\bibitem{KT2}
T.~L. Kriete and D.~Trutt.
\newblock The {C}es\`aro operator in {$l^{2}$} is subnormal.
\newblock {\em Amer. J. Math.}, 93:215--225, 1971.

\bibitem{KT1}
T.~L. Kriete and D.~Trutt.
\newblock On the {C}es\`aro operator.
\newblock {\em Indiana Univ. Math. J.}, 24:197--214, 1974/75.

\bibitem{Lukacs57}
E.~Lukacs.
\newblock Les fonctions caract\'eristiques analytiques.
\newblock {\em Ann. Inst. H. Poincar\'e}, 15:217--251, 1957.

\bibitem{Lukacs}
E.~Lukacs.
\newblock Some extensions of a theorem of {M}arcinkiewicz.
\newblock {\em Pacific J. Math.}, 8:487--501, 1958.

\bibitem{Muntz}
C.~H. M{\"u}ntz.
\newblock {\"U}ber den {A}pproximationssatz con {W}eierstrass.
\newblock In {\em H. A. Schwarz's Festschrift}, pages 303--312. Berlin, 1914.

\bibitem{PW}
R.~E. A.~C. Paley and N.~Wiener.
\newblock {\em Fourier transforms in the complex domain}, volume~19 of {\em
  American Mathematical Society Colloquium Publications}.
\newblock American Mathematical Society, Providence, RI, 1934.
\newblock Reprint of the 1934 original.

\bibitem{PS}
M.~M. Peloso and M.~Salvatori.
\newblock Functions of exponential growth on a half-plane, sets of uniqueness
  and the {M}{\"u}ntz--{S}z{\'a}sz problem for the {B}ergman space.
\newblock {\em preprint}, 2015.

\bibitem{Sed}
A.~M. Sedletskii.
\newblock Complete and incomplete systems of exponentials in spaces with a
  power weight on a half-line.
\newblock {\em Moscow Univ. Math. Bull.}, 69(2):73--76, 2014.
\newblock Translation of Vestnik Moskov. Univ. Ser. I Mat. Mekh. {{\bf{2}}014},
  no. 2, 52--55.

\bibitem{Szasz}
O.~Sz{\'a}sz.
\newblock \"{U}ber die {A}pproximation stetiger {F}unktionen durch lineare
  {A}ggregate von {P}otenzen.
\newblock {\em Math. Ann.}, 77(4):482--496, 1916.

\end{thebibliography}
\bibliographystyle{abbrv}

\end{document}